\documentclass[12pt,a4paper,UTF8]{article}
\usepackage{amsmath,mathtools,amsfonts,amsmath,amssymb,mathrsfs,dsfont,cases,amsthm,bm,bbm}
\usepackage{color,xcolor}
\usepackage{graphicx,subfigure}
\usepackage{titlesec}
\usepackage{geometry}
\usepackage{natbib}
\usepackage[colorlinks=true,citecolor=blue,urlcolor=blue,linkcolor=blue]{hyperref}%
\usepackage{array,booktabs,makecell,tabularx,multirow,diagbox}  

\newcommand{\bR}{\mathbb{R}}
\newcommand{\bP}{\mathbb{P}}

\newcommand{\bbeta}{\bm{\beta}}

\newcommand{\diag}{\mathrm{diag}}

\newtheoremstyle{thry}
{1em}
{1em}
{\itshape\rmfamily}
{}
{\scshape\large}
{.}
{.5em}
{}
\theoremstyle{thry}

\newtheorem{theorem}{\indent Theorem}
\newtheorem{lemma}{\indent  Lemma}

\numberwithin{equation}{section}

\title{Hypothesis testing for homogenous of nodes in $\beta$-models} 
\author{
\textsc{Kang Fu}\thanks{School of Mathematics and Statistics, and Key Laboratory of Nonlinear Analysis \& Applications (Ministry of Education), Central China Normal University, Wuhan 430079, China. \texttt{Email:} \textit{fukang@mails.ccnu.edu.cn}.},
\hspace{1mm}
\textsc{Jianwei Hu}\thanks{School of Mathematics and Statistics, and Hubei Key Laboratory of Mathematical Sciences, Central China Normal University, Wuhan 430079, China. \texttt{Email:} \textit{jwhu@ccnu.edu.cn}.}
\hspace{1mm}
and
\hspace{1mm}
\textsc{Meng Sun}\thanks{School of Mathematics and Statistics, Central China Normal University, Wuhan 430079, China. \texttt{Email:} \textit{sara060@mails.ccnu.edu.cn}.}\\
\textit{Central China Normal University}
}
\date{\empty}

\begin{document}
\maketitle

\begin{abstract}
The $\beta$-model has been extensively utilized to model degree heterogeneity in networks, wherein each node is assigned a unique parameter. In this article, we consider the hypothesis testing problem that two nodes $i$ and $j$ of a $\beta$-model have the same node parameter. We prove that the null distribution of the proposed statistic converges in distribution to the standard normal distribution. Further, we investigate the homogeneous test for $\beta$-model by combining individual $p$-values to aggregate small effects of multiple tests. Both simulation studies and real-world data examples indicate that the proposed method works well.

\end{abstract}

\section{Introduction}

Network models are commonly popular models to character the interaction between the different entries \citep{Scott:2000}. The studies on network data have attracted considerable attention in many fields, such as computer science, social science, and biology. For example, in the social network, the interaction between the different individuals represents a friend relationship \citep{Hunter:2012}. In general, an undirected and unweight network $\mathcal{G}$ with $n$ nodes can be represented by an $n\times n$ adjacency matrix $A\in\{0,1\}^{n\times n}$, where $(i,j)$-th entry indicates whether there is a connection between node $i$ and node $j$, i.e., $A_{ij}=1$ if there is a connection between node $i$ and node $j$ and $A_{ij} = 0$ otherwise. In network data analysis, the $\beta$-model, proposed by \cite{Chatterjee:2011}, is a special case of a class of models known as node-parameter models, where each node degree is associated with a corresponding parameter. Specifically, the $\beta$-model assumes that the edge between node $i$ and node $j$ exists with probability 
\[
\bP\{A_{ij}=1\} = p_{ij} = \dfrac{e^{\beta_i+\beta_j}}{1+e^{\beta_i+\beta_j}},
\]
independently of all other edges, where $\beta_i$ is the node parameter (also known as the ``attractiveness" of node) of node $i$. The $\beta$-model is an exponential random graph model and can be seen as an undirected version of a $p_1$-model \citep{Holland:1981}. An advantage of the $\beta$-model is that the degree sequence is the unique sufficient statistic. Then, the $\beta$-model is widely used to model the network with degree heterogeneous. It is not difficult to see that the probability connecting the node $i$ and node $j$ only depends on the parameters of the node $i$ and node $j$. When all $\beta_i$'s are equal to each other, the $\beta$-model naturally degenerates to the E-R model. To fit a sparse network, \cite{Mukherjee:2018} proposed the adjusted $\beta$-model
\[
p_{ij} = \dfrac{\lambda}{n}\dfrac{e^{\beta_i+\beta_j}}{1+e^{\beta_i+\beta_j}},
\]
where $\lambda\in(1,n)$ is used to measure the sparsity of the graph. Since the $\beta$-model can capture important features of real-world networks, the $\beta$-model and its variations have been studied widely in recent years \citep{Chatterjee:2011,Yan:2013,Rinaldo:2013,Ogawa:2013,Yan:2015,Yan:2016}.

Hypothesis testing plays a critical role in the studies on network data \citep{Fu:2022,Fu:2023}. One significant application is to recover the community structure of a network. \cite{Bickel:2016} and \cite{Dong:2020} used the spectral statistic of the normalized adjacency matrix to test whether the network has a community structure, i.e., $H_0: k=1$ for stochastic block models. Then, \cite{Cammarata:2023} considered the global testing problem under the framework of degree-corrected mixed membership models. Further, a majority of methods of the goodness-of-fit test for stochastic block models have also been proposed, see, e.g., \cite{Lei:2016,Hu:2021,Jin:2023jasa}. Under the settings of degree-corrected mixed membership models, \cite{Fan:2022} studied the issue of hypothesis testing for the equality of membership vectors between two nodes, up to a possible scaling. Similarly, \cite{Du:2023} investigated the equality of latent positions between two nodes. Their methods are based on the Mahalanobis distance between two vectors, which are generalizations of the corresponding results in \cite{Fan:2022}.

Hypothesis testing for $\beta$-modes is a nascent research area. Motivated by the issues of equality of two nodes, we consider the hypothesis testing problem that two node $i$ and node $j$ of a $\beta$-model have the same node parameter. Specifically, we consider the following test:
\begin{equation}\label{eq:test1}
	H_0: \beta_i=\beta_j\quad v.s.\quad H_1:\beta_i\neq\beta_j,
\end{equation}
for any $i,j\in[n]$, where $[n]=\{1,\cdots,n\}$. Further, the other significant problem is the homogeneous test, i.e.,
\begin{equation}\label{eq:test2}
	H_0^{\prime}: \beta_1=\beta_2=\cdots=\beta_n\quad v.s.\quad H_1^{\prime}:\beta_i\neq\beta_j\ \text{for at least one pair of}\  i,j.
\end{equation}
For test \eqref{eq:test2}, the null hypothesis $H_0^{\prime}$ implies that there is no heterogeneity in the network, and the network can be seen as an E-R graph. For an adjusted $\beta$-model, \cite{Mukherjee:2018} considered a homogeneous null hypothesis with all $\beta_i$ being equal to 0 against an alternative hypothesis with a subset of $\{\beta_i:i\in[n]\}$ strictly greater than 0. They proposed three explicitly degree-based test statistics: $\sum_i d_i$, $\max_i d_i$, and a criticism test based on $(d_i-\lambda/2)/(\lambda(1-\lambda/2n))^{1/2}$ and established their asymptotic null distribution under some mild conditions. Similarly, under the $\beta$-model, \cite{Yan:2022} investigated two testing problems: for a fixed $r$, the specified null $H_0:\beta_i=\beta_i^{0}, i=1,\cdots, r$ and the homogeneous null $H_0:\beta_1=\cdots=\beta_r$, where $\beta_i^{0}$ is known. For the two nulls, they established the Wilks' theorem of $\beta$-models, i.e., the log-likelihood ratio statistic $2[\ell(\hat{\bbeta})-\ell(\hat{\bbeta}^{res})]$ converges in distribution to a chi-square distribution with $r$ degrees of freedom and $r-1$ degrees of freedom, respectively, where $\hat{\bbeta}$ and $\hat{\bbeta}^{res}$ are the unrestricted and restricted maximum likelihood estimators of $\bbeta$, and $\ell(\cdot)$ is the log-likelihood function. Compared with their settings, the advantages of our setting are as follows. First, our null hypothesis has a wider range of parameters than that in \cite{Mukherjee:2018} since we do not require that all parameters be equal to zero. Second, we only need the unrestricted maximum likelihood estimate, and save the computational cost.

The rest of this article is organized as follows. In Section \ref{sec:test1}, we present our main method and theorems about the test for equality of node parameters. The homogeneous test for the $\beta$-model is investigated in Section \ref{sec:test2}. Additional simulation studies and real-world data examples are given in Sections \ref{sec:sim} and \ref{sec:real}. Section \ref{sec:conclusion} concludes the article. Technical proofs are given in the Appendix.

\section{Hypothesis testing for equality of node parameters}\label{sec:test1}
Formally, suppose that $A\in\{0,1\}^{n\times n}$ is an adjacency matrix of undirected graph $\mathcal{G}$ generated from the $\beta$-model with parameter $\bbeta = (\beta_1,\ldots,\beta_n)^\top\in\bR^{n}$, where $\bbeta$ is unknown. Throughout this article, we assume that the self-loops are not allowed, i.e., $A_{ii}=0$ for $1\leq i\leq n$. Let $d_i=\sum_{j\neq i}A_{ij}$ be the degree of the node $i$. Then, the logarithm of the likelihood function can be written as:
\begin{equation*}
	\ell(\bbeta|A)=\sum_{i}\beta_id_i-\sum_{i<j}\log\left(1+e^{\beta_i+\beta_j}\right).
\end{equation*}
Denote $\hat{\bbeta}=\mathop{\arg\max}_{\bbeta}\ell(\bbeta|A)$ as the maximum likelihood estimator (MLE). The MLE can be obtained by solving the following equations:
\begin{equation}\label{eq:MLE}
	d_i=\sum_{j\neq i}\dfrac{e^{\beta_i+\beta_j}}{1+e^{\beta_i+\beta_j}}, (i=1,\cdots,n).
\end{equation}
\cite{Chatterjee:2011} showed that the fixed point iterative algorithm can be used to solve $\hat{\bbeta}$. Under the frameworks of the $\beta$-model, \cite{Chatterjee:2011} established the consistency of $\hat{\bbeta}$. Specifically, let $L_n=\max_{i}|\beta_i|$, then there is a constant $C(L_n)$ depending only on $L_n$ such that $\bP\{\max_{1\leq i\leq n}|\hat{\beta}_i-\beta_i|\leq C(L_n)\sqrt{n^{-1}\log n}\}\geq 1-C(L_n)n^{-2}$. Further, by approximating the inverse of the Fisher information matrix, \cite{Yan:2013} proved the asymptotic normality of $\hat{\bbeta}$. Then, \cite{Rinaldo:2013} gave the necessary and sufficient conditions for the existence and uniqueness of $\hat{\bbeta}$.

Denote the Fisher information matrix for $\bbeta$ as $V=(v_{ij})_{n\times n}$, where 
\[
v_{ij}=\dfrac{e^{\beta_i+\beta_j}}{\{1+e^{\beta_i+\beta_j}\}^2}\ (1\leq i\neq j\leq j), \qquad v_{ii}=\sum_{j\neq i}v_{ij}.
\]
Note that $V$ is also the covariance matrix of degree sequence $\bm{d}=(d_1,\cdots,d_n)$. Then, \cite{Yan:2013} established the following central limiting theorem:
\begin{lemma}\label{thm:yan}
	If $L_n=o(\log\log n)$, then for any fixed $r\geq 1$, the vector consisting of the first $r$ elements of $G^{1/2}(\hat{\bbeta}-\bbeta)$ is asymptotically standard multivariate normal as $n\rightarrow\infty$, where $G=\diag(v_{11},\cdots ,v_{nn})$ and $G^{1/2}=\diag(v_{11}^{1/2},\cdots ,v_{nn}^{1/2})$.
\end{lemma}

Lemma \ref{thm:yan} implies that, for any $i\in[n]$, the following result holds:
\[
v_{ii}^{1/2}(\hat{\beta}_i-\beta_i)\stackrel{d}{\longrightarrow} N(0,1),
\]
and $\hat{\beta}_i$ and $\hat{\beta}_j$ are asymptotic independent for any $1\leq i\neq j\leq n$. Then, for a pair of nodes $(i,j)$,  we have
\[
\hat{\beta}_i-\hat{\beta}_j\stackrel{d}{\longrightarrow} N(\beta_i-\beta_j,v_{ii}^{-1}+v_{jj}^{-1}).
\]
Under the null hypothesis $H_0$ of test \eqref{eq:test1}, we have
\begin{equation}\label{eq:betec}
	\hat{\beta}_i-\hat{\beta}_j\stackrel{d}{\longrightarrow} N(0,v_{ii}^{-1}+v_{jj}^{-1}).
\end{equation}
Consider the statistic $U_{ij}=\dfrac{\hat{\beta}_i-\hat{\beta}_j}{\sqrt{v_{ii}^{-1}+v_{jj}^{-1}}}$. Then, under $H_0$, we have $U_{ij}\stackrel{d}{\longrightarrow} N(0,1)$. Notice that the statistic $U_{ij}$ involves unknown parameters $v_{ii}$ and $v_{jj}$. Hence, we can consider a natural estimate of $U_{ij}$ by plugging in the estimated parameters $\hat{v}_{ii}$ and $\hat{v}_{jj}$, where
\[
\hat{v}_{ij}=\dfrac{e^{\hat{\beta}_i+\hat{\beta}_j}}{\{1+e^{\hat{\beta}_i+\hat{\beta}_j}\}^2}\ (1\leq i\neq j\leq n), \qquad \hat{v}_{ii}=\sum_{j\neq i}\hat{v}_{ij}.
\]
Denote the empirical estimate of $U_{ij}$ by $\hat{U}_{ij}=\dfrac{\hat{\beta}_i-\hat{\beta}_j}{\sqrt{\hat{v}_{ii}^{-1}+\hat{v}_{jj}^{-1}}}$. It is natural to conjecture that when the estimates $\hat{v}_{ii}$ and $\hat{v}_{jj}$ are accurate enough, the convergence in \eqref{eq:betec} will still hold for $\hat{U}_{ij}$.

Formally, we have the following theorem:
\begin{theorem}\label{thm:main}
	Let $A$ be an adjacency matrix generated from a $\beta$-model with parameter $\bbeta=(\beta_1,\cdots,\beta_n)$. Under $H_0$, when $\max_{i}|\beta_i|=o(\log\log n)$, we have the following result:
	\begin{equation}
		\hat{U}_{ij}\stackrel{d}{\longrightarrow} N(0,1).
	\end{equation}
	
	Under $H_1$, we assume that $\beta_i-\beta_j=\mu$. Then, we have
	\begin{equation}
		\hat{U}_{ij}\stackrel{d}{\longrightarrow} N(\mu,1).
	\end{equation}
\end{theorem}
We postpone the proof to the Appendix. Theorem \ref{thm:main} is an intuitive result. The method is similar to the test of the mean for two samples when the variance is unknown. It can be seen that, for the null and alternative, the statistic $\hat{U}_{ij}$ has different means. Using the result, we can carry out the hypothesis testing. Specifically, given a nominal level $\alpha$, we have a rejection rule:
\begin{equation}
	\mathrm{Reject}\ H_0,\ \mathrm{if}\ |\hat{U}_{ij}|\geq u_{1-\alpha/2},
\end{equation}
where $u_{1-\alpha/2}$ is the upper $\alpha$-th quantile of the standard normal distribution.

\section{Hypothesis testing for homogeneous}\label{sec:test2}

In this section, we consider the homogeneous testing for the $\beta$-model. Under the null hypothesis of test \eqref{eq:test2}, the $\beta$-model reduces to the E-R model. Then, the homogeneous testing enables the evaluation of heterogeneity among the nodes within the network. For the test \eqref{eq:test2}, the alternative hypothesis implies that there is a pair of nodes $(i,j)$ with non-equality of node parameters. Hence, using the test \eqref{eq:test2} on node pairs $(i, j)$ will result in rejecting the null hypothesis. Intuitively, we can consider all pairs of nodes $(i,j)$ for $1\leq i < j\leq n$, then using the test \eqref{eq:test2} on node pairs $(i, j)$, which leads to $n(n-1)/2$ testing results. A significant problem is the statistics $U_{ij}$'s are correlated and how to combine the information of $n(n-1)/2$ results.

In the meta-analysis, methods for combining multiple test statistics are widely used in massive data analysis. Specifically, suppose we independently test the same hypothesis using $K$ different statistical tests and obtain $p$-values $p_1,\cdots, p_K$. An important issue is how to combine them into a single $p$-value. Notice that, under the null hypothesis, all $p_i$'s should follow the uniform distribution on interval $[0,1]$. Hence, the null hypothesis can be rewritten as 
\[
H^{''}_0: p_i\sim U[0,1]\ \text{for}\ i=1,\cdots,K.
\]
The six most simple and commonly used statistics for combining $p$-values are: $T_F=\sum_{i}\log p_i$ \citep{Fisher:1932}, $T_P=-\sum_{i}\log(1-p_i)$ \citep{Pearson:1993}, $T_G=T_F+T_P=\sum_i\log\{p_i/(1-p_i)\}$ \citep{Mudholkar:1979}, $T_E=\sum_i p_i$ \citep{Edgington:1972}, $T_S=\sum_i\Phi^{-1}(p_i)$ \citep{Stouffer:1949}, $T_T=\min_i p_i$ \citep{Tippett:1931}. However, an obvious deficiency is that, when there is a dependence structure between $p_i$'s, all these six methods do not work. Then, \cite{Liu:2020} proposed a Cauchy combination method that takes advantage of the Cauchy distribution. A nonasymptotic result was established to demonstrate that the tail of the null distribution can be effectively approximated by a Cauchy distribution, under arbitrary dependency structures. Specifically, the Cauchy test statistic has the form: $T_L = \sum_i w_i\tan\{(0.5-p_i)\cdot\pi\}$, where the weights $w_i$'s are nonnegative and $\sum_i w_i=1$.

Recall the homogeneous test \eqref{eq:test2}. For any pair of nodes $(i,j)$, we can calculate the statistic $\hat{U}_{ij}$ and the $p$-value $p_{ij}^{value}=2\bP_{N(0,1)}\{X\geq|\hat{U}_{ij}|\}$. Under the null $H_{0}^{\prime}$, all $p_{ij}^{value}$'s should follow the uniform distribution on interval $[0,1]$, and they are not independent. Hence, we consider the Cauchy combination statistic:
\[
T_n=\sum_{i< j}w_{ij}\tan\{(0.5-p_{ij}^{value})\cdot\pi\}.
\]
According to the results in \cite{Liu:2020}, the test statistic $T_n$ has approximately a Cauchy tail even when $p_{ij}^{value}$'s are dependent, i.e.,
\[
\lim\limits_{t\rightarrow+\infty}\dfrac{\bP\{|T_n|\geq t\}}{\bP\{|C_0|\geq t\}} =1,
\]
where $C_0$ denotes a standard Cauchy random variable. Then, for a given nominal level $\alpha$, we have the reject rule:
\[
\text{Reject}\ H_0^{\prime}:\ \text{if}\ |T_n|>c_{1-\alpha/2},
\]
where $c_{1-\alpha/2}$ is the upper $\alpha$-th quantile of the standard Cauchy distribution.

\textbf{Remark.} Compared with the resluts in \cite{Yan:2022}, the proposed method can test the homogeneous for $n$ parameters. Following Lemma \ref{thm:yan}, when $r$ diverges to infinity, the first $r$ elements of $\hat{\bbeta}$ may not be independent. However, in our test procedure, we only consider the two estimators $\hat{\beta}_i$ and $\hat{\beta}_j$ that can be seen as independent, then we can combine information from $n(n-1)/2$ tests.

\section{Simulation}\label{sec:sim}
In this section, we carry out extensive simulation studies to evaluate the performance of the proposed method. All simulations were performed on a PC with a single processor of 2.3 GHz 8‐Core Intel Core i9.

\subsection{The empirical distribution for statistic $\hat{U}_{ij}$}\label{sim:1}
In this simulation, we examine the finite sample empirical distribution of the test statistic $\hat{U}_{ij}$ under the null and alternative hypothesis and verify the result in Theorem \ref{thm:main}. We set $n=300, 500$ and $\beta_i=iL_n/(n-1)$ where $L_n=0, \log(\log n)$, and $(\log n)^{1/2}$. When $L_n=0$, all $\beta_i$'s are equal, which corresponds to $H_0$. And, when $L_n=\log(\log n)$ and $(\log n)^{1/2}$, there is heterogeneous between nodes, which corresponds to $H_1$.

\begin{figure}[htbp]
	\centering
	\includegraphics[width=\textwidth]{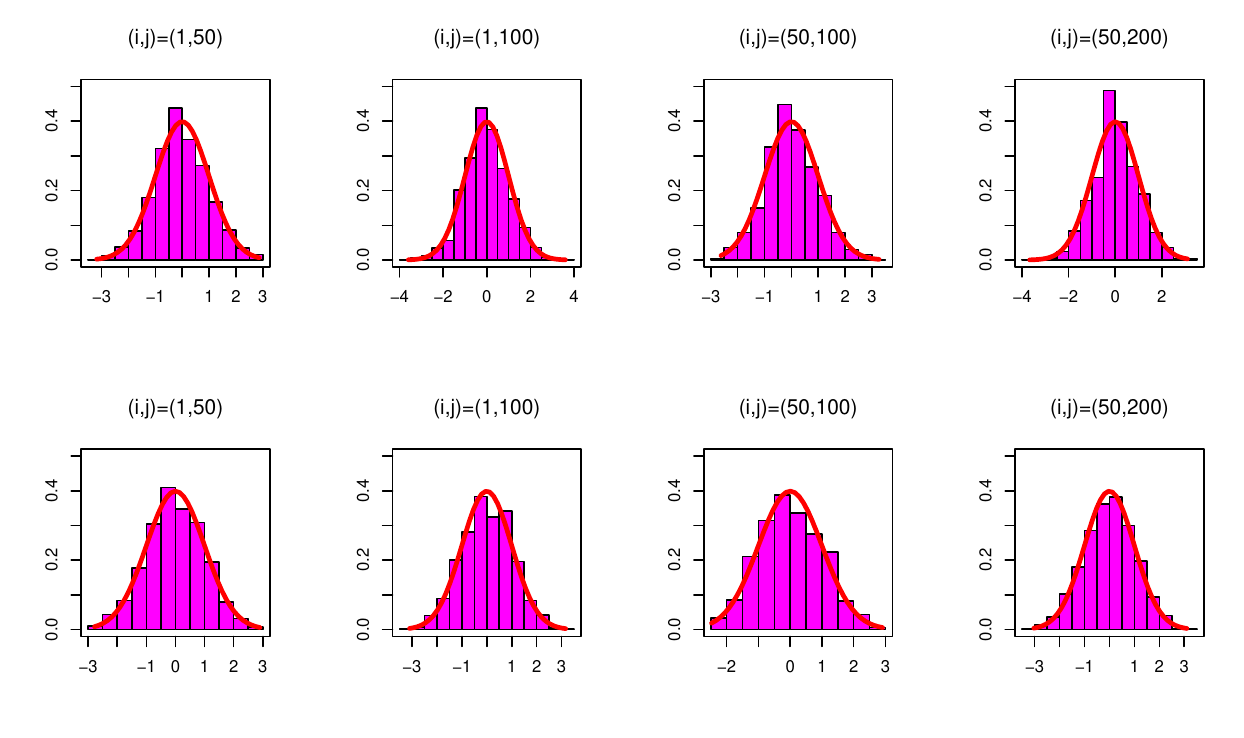}
	\caption{The histogram of the statistic $\hat{U}_{ij}$ under $n=300$ (upper row) and $n=500$ (lower row) when $L_n = 0$. The red solid line indicates the density of the standard normal distribution.}
	\label{fig:density1}
\end{figure}

\begin{figure}[htbp]
	\centering
	\includegraphics[width=\textwidth]{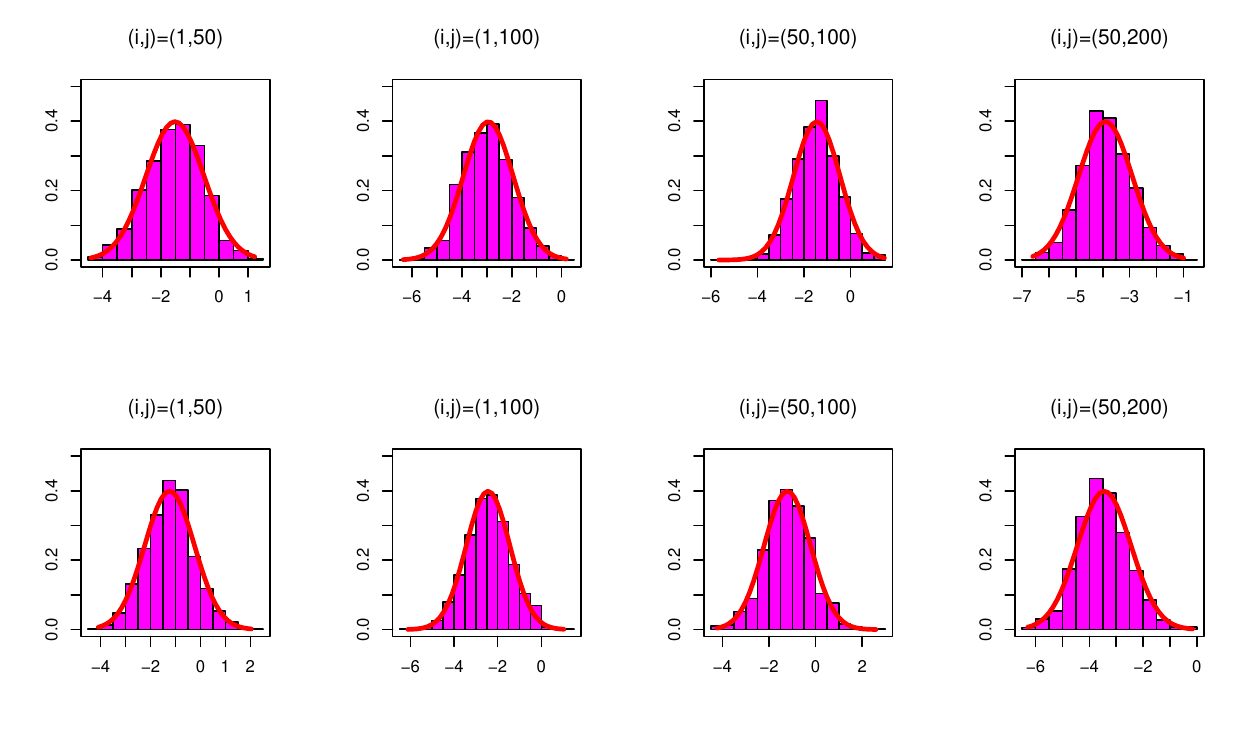}
	\caption{The histogram of the statistic $\hat{U}_{ij}$ under $n=300$ (upper row) and $n=500$ (lower row) when $L_n = \log\log(n)$. The red solid line indicates the density of the normal distribution with $\mu=\beta_i-\beta_j$ and $\sigma^2=1$.}
	\label{fig:density2}
\end{figure}

\begin{figure}[htbp]
	\centering
	\includegraphics[width=\textwidth]{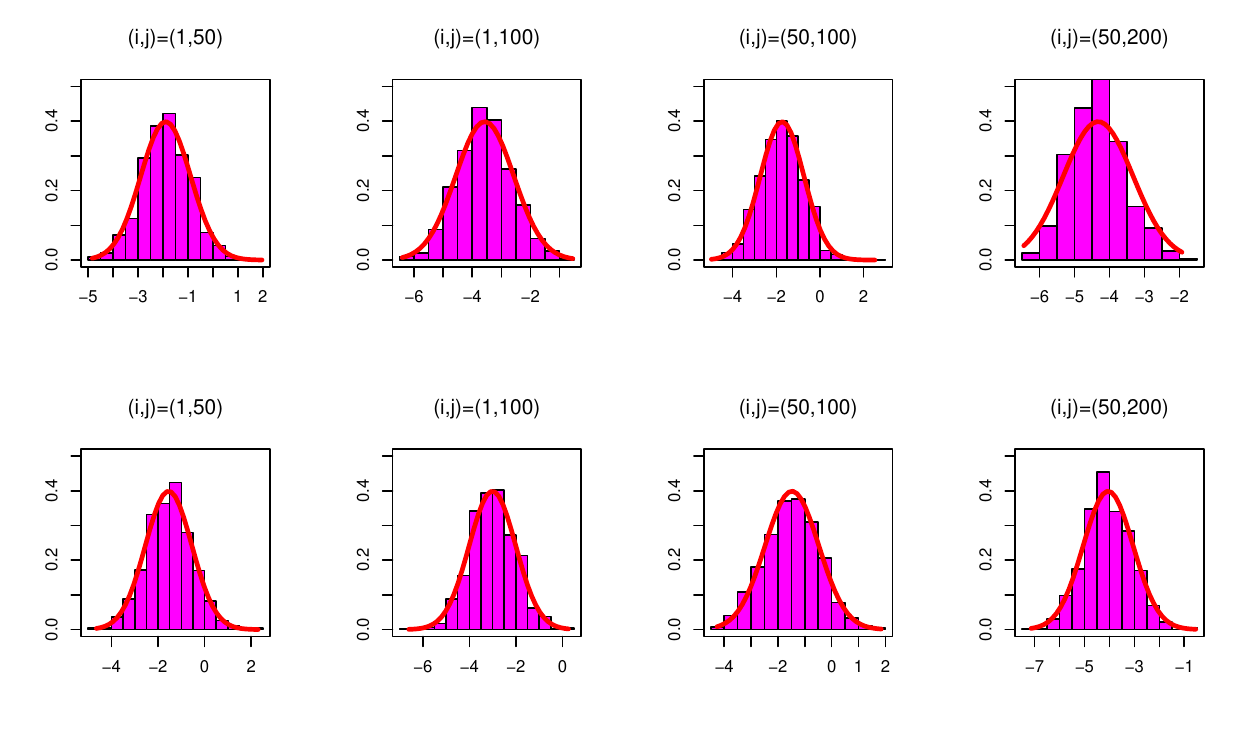}
	\caption{The histogram of the statistic $\hat{U}_{ij}$ under $n=300$ (upper row) and $n=500$ (lower row) when $L_n = (\log(n))^{1/2}$. The red solid line indicates the density of the standard normal distribution with $\mu=\beta_i-\beta_j$ and $\sigma^2=1$.}
	\label{fig:density3}
\end{figure}

In Figures \ref{fig:density1}-\ref{fig:density3}, we plot the empirical density of the statistic $\hat{U}_{ij}$ from 1000 data replications. When $L_n=0,\log(\log n)$, and $(\log n)^{1/2}$, the plots show that the simulation result very well matches the prediction of Theorem \ref{thm:main}. Under the null ($L_n=0$) and the alternative ($L_n=\log(\log n)$ and $(\log n)^{1/2}$), the test statistic has different mean.

\subsection{The empirical size and power for test \eqref{eq:test1}}
In this subsection, we investigate the empirical size and power for test \eqref{eq:test1}, and the settings are similar to that in Section \ref{sim:1}. The proportion of rejection at nominal level 0.05 is summarized in Table \ref{tab:test1}. It is easy to see that the type I error is correctly kept at the nominal level. For the alternative hypothesis, the power tends to be less than 1. In fact, when the difference between $\beta_i$ and $\beta_j$ is small ($(i,j)=(1,50)$ or $(50,100)$), the distribution of $\hat{U}_{ij}$ is close to the standard normal distribution, which leads to the power may be much less than 1. When the difference between $\beta_i$ and $\beta_j$ is large ($(i,j)=(1,100)$ or $(50,200)$), however, the empirical powers are close to 1. The results are consistent with the results in Section \ref{sim:1}. In addition, we observe that, with the sample increasing, the power of the test decreases. The main reason is that the parameter generation method makes the difference between two nodes become smaller as the number of samples increases.

\begin{table}[htbp]
\setlength{\abovecaptionskip}{0cm}  
\setlength{\belowcaptionskip}{0.5cm} 
\centering
\caption{The proportion of rejection at nominal level 0.05 over 200 independent samples.}
\label{tab:test1}
\begin{tabular}{cccccccc}
\toprule
 & \multicolumn{3}{c}{$n=300$} & & \multicolumn{3}{c}{$n=500$} \\ \cline{2-4}\cline{6-8}
$(i,j)$ & $L_n=0$ & $L_n=\log\log n$ & $(\log n)^{1/2}$ & & $L_n=0$ & $L_n=\log\log n$ & $L_n=(\log n)^{1/2}$ \\ \midrule
$(1,50)$ & 0.05 & 0.34 & 0.47 & & 0.05 & 0.25 & 0.35 \\
$(1,100)$ & 0.05 & 0.86 & 0.96 & & 0.05 & 0.73 & 0.87 \\
$(50,100)$ & 0.06 & 0.31 & 0.41 & & 0.06 & 0.27 & 0.31 \\
$(50,200)$ & 0.07 & 0.95 & 0.99 & & 0.05 & 0.96 & 0.99 \\ \bottomrule
\end{tabular}
\end{table}

\subsection{The empirical size and power for test \eqref{eq:test2}}
In this subsection, we investigate the homogeneous test for the $\beta$-model. We also set $\beta_i=(i-1)L_n/(n-1)$. However, we set $\beta_1=\cdots=\beta_r = 0$, where $r$ has five cases: $n, n-1, n-2, n-5$, and $n-10$. It is easy to see that $r=n$ corresponds to the null $H_0^{\prime}$, and the other four cases correspond to the null $H_1^{\prime}$. For $L_n$, we consider two classes of settings: (i) $L_n=(\log(\log n))^{1/2}, \log(\log n)$, and $(\log n)^{1/2}$; (2) $L_n=c\log n$, where $c=0.1,0.2$, and $0.5$. The results are given in Tables \ref{tab:homo1} and \ref{tab:homo2}. For the simulation results, under the null $H_0^{\prime}$ ($r=n$), the type I errors are close to the nominal level. For the alternative $H_1^{\prime}$, the empirical power is less than 1, and the proposed method is superior to the method in \cite{Yan:2022} when $r$ approximates $n$. All simulation results show that the proposed method is effective and efficient.

\begin{table}[htbp]
\setlength{\abovecaptionskip}{0cm}  
\setlength{\belowcaptionskip}{0.5cm} 
\centering
\caption{The proportion of rejection at nominal level 0.05 over 200 independent samples.}
\label{tab:homo1}
\begin{tabular}{c|l|ccc}
\toprule
&& $L_n=(\log(\log n))^{1/2}$ & $L_n=\log(\log n)$ & $L_n=(\log n)^{1/2}$ \\ \midrule
\multirow{5}{*}{$n=100$}
& $r=n$ & 0.05 (0.07) & 0.03 (0.08) & 0.05 (0.05) \\
& $r=n-1$ & 0.97 (0.54) & 1 (0.81) & 1 (0.99) \\
& $r=n-2$ & 1 (0.96) & 1 (1) & 1 (1) \\
& $r=n-5$ & 1 (1) & 1 (1) & 1 (1) \\ 
& $r=n-10$ & 1 (1) & 1 (1) & 1 (1) \\ \midrule
\multirow{5}{*}{$n=200$}
& $r=n$ & 0.04 (0.08) & 0.03 (0.05) & 0.04 (0.09) \\
& $r=n-1$ & 0.94 (0.52) & 1 (0.82) & 1 (0.99) \\
& $r=n-2$ & 1 (0.94) & 1 (0.99) & 1 (1) \\
& $r=n-5$ & 1 (1) & 1 (1) & 1 (1) \\ 
& $r=n-10$ & 1 (1) & 1 (1) & 1 (1) \\ \midrule
\multirow{5}{*}{$n=500$}
& $r=n$ & 0.01 (0.05) & 0.03 (0.07) & 0.03 (0.05) \\
& $r=n-1$ & 0.98 (0.60) & 1 (0.82) & 1 (0.98) \\
& $r=n-2$ & 0.99 (0.96) & 1 (1) & 1 (1) \\
& $r=n-5$ & 1 (1) & 1 (1) & 1 (1) \\ 
& $r=n-10$ & 1 (1) & 1 (1) & 1 (1) \\
\bottomrule
\end{tabular}
\end{table}

\begin{table}[htbp]
\setlength{\abovecaptionskip}{0cm}  
\setlength{\belowcaptionskip}{0.5cm} 
\centering
\caption{The proportion of rejection at nominal level 0.05 over 200 independent samples.}
\label{tab:homo2}
\begin{tabular}{c|l|ccc}
\toprule
&& $c=0.1$ & $c=0.2$ & $c=0.5$ \\ \midrule
\multirow{5}{*}{$n=100$}
& $r=n$ & 0.05 (0.07) & 0.03 (0.08) & 0.05 (0.05) \\
& $r=n-1$ & 0.09 (0.06) & 0.73 (0.23) & 1 (1) \\
& $r=n-2$ & 0.17 (0.10) & 0.92 (0.67) & 1 (1) \\
& $r=n-5$ & 0.34 (0.35) & 1 (1) & 1 (1) \\ 
& $r=n-10$ & 0.57 (0.78) & 1 (1) & 1 (1) \\ \midrule
\multirow{5}{*}{$n=200$}
& $r=n$ & 0.04 (0.05) & 0.04 (0.06) & 0.05 (0.07) \\
& $r=n-1$ & 0.46 (0.11) & 1 (0.67) & 1 (1) \\
& $r=n-2$ & 0.67 (0.28) & 1 (0.98) & 1 (1) \\
& $r=n-5$ & 0.95 (0.83) & 1 (1) & 1 (1) \\ 
& $r=n-10$ & 1 (1) & 1 (1) & 1 (1) \\ \midrule
\multirow{5}{*}{$n=500$}
& $r=n$ & 0.04 (0.05) & 0.02 (0.06) & 0.03 (0.06) \\
& $r=n-1$ & 0.99 (0.33) & 1 (0.99) & 1 (1) \\
& $r=n-2$ & 1 (0.84) & 1 (1) & 1 (1) \\
& $r=n-5$ & 1 (1) & 1 (1) & 1 (1) \\ 
& $r=n-10$ & 1 (1) & 1 (1) & 1 (1) \\ 
\bottomrule
\end{tabular}
\end{table}

\section{Real example analysis}\label{sec:real}
In this section, we apply the proposed method to a real network dataset. The food web dataset is from \cite{Baird:1989} and is available in \cite{Blitzstein:2011}, which contains data on 33 organisms (such as bacteria, oysters, and catfish) in the Chesapeake Bay during the summer. The degree sequence of this network is $\bm{d} = (7, 8, 5, 1, 1, 2, 8, 10, 4, 2, 4, 5, 3, 6, 7, 3, 2, 7, 6, 1, 2, 9, 6, 1, 3, 4, 6, 3, 3, 3, 2, 4, 4)$. We observe that some nodes have identical degrees in this network, and the heterogeneity of the network seems not very obvious. To investigate the equality of node parameters, we consider the nodes 4, 6, 13, 11, 12, 14, 15, 2, 22, and 8, which correspond to degrees 1, 2, 3, 4, 5, 6, 7, 8, 9, and 10. Table \ref{tab:real} shows that the $p$-values for test problem \eqref{eq:test1}. The result indicates that the increase in degree difference between two nodes leads to a decrease in $p$-value, which tends to reject the null hypothesis. Finally, we consider the homogeneous test \eqref{eq:test2}. The $p$-values obtained by the proposed method and likelihood-ratio test are 0.698 and 0.998, respectively. The result shows that the network is homogeneous with high probability. 

\begin{table}[htbp]
\setlength{\abovecaptionskip}{0cm}  
\setlength{\belowcaptionskip}{0.5cm} 
\centering
\caption{The $p$-values of the test statistic $\hat{U}_{ij}$ under the test problem \eqref{eq:test1}.}
\label{tab:real}
\begin{tabular}{c|cccccccccc}
\diagbox{$i$}{$j$} & 4 & 6 & 13 & 11 & 12 & 14 & 15 & 2 & 22 & 8 \\ \hline
4 & $-$ & 0.277 & 0.156 & 0.090 & 0.053 & 0.031 & 0.019 & 0.012 & 0.007 & 0.004 \\
6 & 0.277 & $-$ & 0.316 & 0.189 & 0.110 & 0.063 & 0.035 & 0.019 & 0.011 & 0.006 \\
13 & 0.156 & 0.316 & $-$ & 0.337 & 0.213 & 0.128 & 0.074 & 0.042 & 0.023 & 0.012 \\
11 & 0.090 & 0.189 & 0.337 & $-$ & 0.350 & 0.230 & 0.143 & 0.085 & 0.049 & 0.027 \\
12 & 0.053 & 0.110 & 0.213 & 0.350 & $-$ & 0.360 & 0.243 & 0.156 & 0.095 & 0.055 \\
14 & 0.031 & 0.063 & 0.128 & 0.230 & 0.360 & $-$ & 0.367 & 0.254 & 0.167 & 0.104 \\
15 & 0.019 & 0.035 & 0.074 & 0.143 & 0.243 & 0.367 & $-$ & 0.373 & 0.263 & 0.176 \\
2  & 0.012 & 0.019 & 0.042 & 0.085 & 0.156 & 0.254 & 0.373 & $-$ & 0.378 & 0.271 \\
22 & 0.007 & 0.011 & 0.023 & 0.049 & 0.095 & 0.167 & 0.263 & 0.378 & $-$ & 0.382 \\
8  & 0.004 & 0.006 & 0.012 & 0.027 & 0.055 & 0.104 & 0.176 & 0.271 & 0.382 & $-$ \\ \hline
\end{tabular}
\end{table}

\section{Conclusion}\label{sec:conclusion}

In this article, we have proposed a novel statistic to investigate the equality test for the two nodes of the $\beta$-model. Based on the central limit theorem, we have proved the limiting distribution of the proposed statistic is the standard normal distribution. Then, plugging in the MLE of parameters, we have proved that the limiting distribution of the empirical counterpart of the test statistic is also the standard normal distribution under some mild conditions. Under the alternative hypothesis, the limit distribution of the test statistic has also been proven to be a normal distribution with a different mean from the null distribution. Further, based on the combining $p$-values method, we have investigated the homogeneous test for the $\beta$-model. Empirically, by extensive simulation studies, we have demonstrated that the size and the power of the test are valid.

It is worth noting that the proposed test method works well when the difference between the parameters of two nodes is large. However, the power will decrease when the difference between the parameters of two nodes is small. Hence, we need to consider how to improve the power of the proposed test for hypothesis test \eqref{eq:test1} under the case of $0<\beta_i-\beta_j\leq \varepsilon$ for a small constant $\varepsilon>0$. Next, we can also consider extending the single sample to the multi-sample, such as $H_0: \bbeta_1 = \bbeta_2$ for two $\beta$-models with parameters $\bbeta_1$ and $\bbeta_2$. We will continue to study this issue in future work.

\section{Appendix}

\subsection{Proof of Theorem \ref{thm:main}}

First, we consider the case of $H_0$. According to the Taylor expansion, we have, for any $1\leq i\leq n$,
\[
\hat{v}_{ii}^{-1}-v_{ii}^{-1} = -v_{ii}^{-2}(\hat{v}_{ii}-v_{ii}).
\]
Following the definition of $v_{ij}$, it is easy to see that
\begin{equation}\label{eq:proof1}
	\dfrac{n-1}{4}e^{-2L_n}\leq v_{ii}\leq \dfrac{n-1}{4}\quad \text{and}\quad \dfrac{16}{(n-1)^2}\leq v_{ii}^{-2}\leq \dfrac{16}{(n-1)^2}e^{-4L_n}.
\end{equation}

Next, we consider to bound the terms $\hat{v}_{ii}-v_{ii}$. Define $f(x)=e^x/\{1+e^x\}^2$, then $f^{\prime}(x)=-e^x(e^x-1)/\{1+e^x\}^3$. For any $1\leq i\leq n$,
\begin{align*}
	\hat{v}_{ii}-v_{ii} & \leq \sum_{j\neq i}|\hat{v}_{ij}-v_{ij}| \\
	& = \sum_{j\neq i}|f(\hat{\beta}_i+\hat{\beta}_j)-f(\beta_i+\beta_j)| \\
	& = \sum_{j\neq i}|f^{\prime}(\beta_i+\beta_j)(\hat{\beta}_i+\hat{\beta}_j-\beta_i-\beta_j)| \\
	& \leq \sum_{j\neq i}2|f^{\prime}(\beta_i+\beta_j)|\cdot|\hat{\beta}_i-\beta_i|.
\end{align*}
Notice that $|f^{\prime}(x)|\leq 1/6\sqrt{3}$ and the convergence rate of $\hat{\beta}_i$ is between $O_p(n^{-1/2}e^{L_n})$ and $O_p(n^{-1/2})$. Hence, we have $|\hat{v}_{ii}-v_{ii}|=O_p(n^{1/2}e^{L_n})$. Combining with \eqref{eq:proof1}, we have, for any $1\leq i\leq n$, 
\[
|\hat{v}_{ii}^{-1}-v_{ii}^{-1}| = O_p(n^{-3/2}e^{5L_n}).
\]

Thus, we have, for any $1\leq i\neq j\leq j$,
\begin{align*}
	\hat{U}_{ij} & =\dfrac{\hat{\beta}_i-\hat{\beta}_j}{\sqrt{\hat{v}_{ii}^{-1}+\hat{v}_{jj}^{-1}}}\\
	& = \dfrac{\hat{\beta}_i-\hat{\beta}_j}{\sqrt{v_{ii}^{-1}+v_{jj}^{-1}}}\times\dfrac{\sqrt{v_{ii}^{-1}+v_{jj}^{-1}}}{\sqrt{\hat{v}_{ii}^{-1}+\hat{v}_{jj}^{-1}}} \\
	& = \dfrac{\hat{\beta}_i-\hat{\beta}_j}{\sqrt{v_{ii}^{-1}+v_{jj}^{-1}}}\times\left(1+O_p(n^{-3/4}e^{5L_n/2})\right).
\end{align*}
According to the Slutsky's theorem, we have $\hat{U}_{ij}\stackrel{d}{\longrightarrow}N(0,1)$.

The proof of the alternative $H_1$ are similar to that of the null $H_0$, we omit the details in the article.

\section*{Acknowledgments}
Hu is partially supported by the National Natural Science Foundation of China (nos. 12171187, 12371261).

\bibliographystyle{agu04}
\bibliography{ref}

\end{document}